\documentclass{article}

\usepackage{amsmath,amssymb}
\usepackage{graphicx}
\usepackage{color}
\usepackage{dsfont}

\newtheorem{theorem}{Theorem}[section]%
\newtheorem{corollary}[theorem]{Corollary}%
\newtheorem{proposition}[theorem]{Proposition}%

\newtheorem{definition}[theorem]{Definition}%
\newtheorem{remark}[theorem]{Remark}%
\newtheorem{example}[theorem]{Example}%

\newcommand{\N}{\ensuremath{\mathbb N}}

\newcommand{\done}{\hfill $\Box$ }

\newcommand{\ls}[1]
    {\dimen0=\fontdimen6\the\font\lineskip=#1\dimen0
     \advance\lineskip.5\fontdimen5\the\font
     \advance\lineskip-\dimen0
     \lineskiplimit=0.9\lineskip
     \baselineskip=\lineskip
     \advance\baselineskip\dimen0
     \normallineskip\lineskip\normallineskiplimit\lineskiplimit
     \normalbaselineskip\baselineskip
     \ignorespaces}

\begin{document}

\bibliographystyle{abbrv}

\title{The $L^0$-extension of an $L^\infty$-normed module}
\author{Mingzhi Wu \quad Tiexin Guo$^\ast$\\
School of Mathematics and Statistics \\
Central South University\\
Changsha {\rm 410083}, China \\
Email: wumingzhi@csu.edu.cn, tiexinguo@csu.edu.cn\\
}

\date{}
 \maketitle

\thispagestyle{plain}
\setcounter{page}{1}

\begin{abstract}
In this paper, we embed each $L^\infty$-normed module $E$ into an appropriate and unique complete random normed module $E_0$ so that the properties of $E$ are closely related to the properties of $E_0$.


\end{abstract}

\ls{1.5}
\section{Introduction}
Since Artzner et al. presented
and studied coherent risk measures in the seminal paper [1], the theory of risk measures has obtained a quite extensive development in the last 15 years: from coherent risk measures \cite{ADEH,Delbaen} to convex risk measures \cite{FS1,FS2,FR1,FR2} and to conditional and dynamic risk measures \cite{Bion-Nadal,Det-Scan,FKVA} and so on. For some most basic concepts and results, see Chapter 4 of the textbook \cite{FS}.

In the static case, the utmost general model space to define a risk measure is a topological vector spaces, functional analysis and convex analysis play an important role in the analysis of risk measures in such kind. In the conditional case, topological vector spaces are no longer suitable model spaces to define a conditional risk measure. Filipovi\'c et al \cite{FKV,FKVA} have shown that it is natural to study conditional risk measures in module framework, where the module is over the ring $L^0$ of all random variables. The module approach to conditional risk measures needs to extend known results in functional analysis and convex analysis from topological vector spaces to topological $L^0$-modules.

The most important topological $L^0$-modules are random normed modules and random locally convex modules, the study of which went to back to the early 1990s by Guo. In the past 20 plus years, the theory of random normed modules and random locally convex modules have undergone a systematic and deep
development: the most basic notions are formulated and lots of important results are obtained. We recommend \cite{Guo-JFA,Guo-recent} on which the most
important results are resumed. Especially, in \cite{Guo-JFA}, Guo gave the relations between some basic results relevant for the financial applications for the random locally convex modules-Guo's results under the $(\varepsilon,\lambda)$-topology and Filipovi\'c et al's results under the locally $L^0$-convex topology.  Now, by considering the two kinds of topologies simultaneously, making full use of the advantage of each kind of topology, both the theory of random locally convex modules and its financial application are developing rapidly.

Recently, Eisele and Taieb \cite{ET1,ET2} made an endeavor to module approach in a somewhat different direction by choosing modules with $L^\infty$ spaces as ring. They succeed in extending many important results in functional analysis to topological $L^\infty$-modules.

In this paper, we try to connect Eisele and Taieb's study for $L^\infty$-modules and Guo's study for random normed modules. We start by embedding a given $L^\infty$-normed module $E$ into an appropriate random normed module $E_0$ in Section 3, then in Section 4, we show that some characterizations for the reflexivity and subreflexivity of $E$ can be derived from the random reflexivity and random subreflexivity of $E_0$.

\section{Terminology and notation}

Let $(\Omega,{\mathcal F},P)$ be a probability space, $L^{0}(\mathcal{F})$ (${\bar L}^{0}(\mathcal{F})$) be the algebra of all equivalence classes of ${\mathcal F}$--measurable real valued (accordingly, extended real valued) random variables on $\Omega$, and $L^\infty(\mathcal{F})$ be the algebra of all equivalence classes of essentially bounded real valued random variables. If there is no other $\sigma$-algebra to be considered, we write $L^0$ and $L^\infty$ for $L^{0}(\mathcal{F})$ and $L^\infty(\mathcal{F})$, respectively.

As usual, $L^{0}$ is partially ordered by $\xi\leqslant\eta$ iff $\xi^{0}(\omega)\leq\eta^{0}(\omega)$ for $P$--almost all $\omega\in \Omega$, where $\xi^0$ and $\eta^0$ are arbitrarily chosen representatives of $\xi$ and $\eta$, respectively. According to \cite{Dunford}, $(L^0,\leqslant)$ is a conditionally complete lattice. For a subset $A$ of $L^0$ with an upper bound (a lower bound), $\vee A$ (accordingly, $\wedge A$) stands for the supremum (accordingly, infimum) of $A$.

${\tilde I}_A$ always denotes the equivalence class of $I_A$, where $A\in {\mathcal F}$ and $I_A$ is the characteristic function of $A$. For any $\xi\in L^0$, $|\xi|$ denotes the equivalence class of $|\xi^0|: \Omega\to [0,\infty)$ defined by $|\xi^0|(\omega)=|\xi^0(\omega)|$, where $\xi^0$ is an arbitrarily chosen representative of $\xi$.

Denote $L^{0}_{+}=\{\xi\in L^{0}\,|\,\xi\geqslant 0\}$ and $L^{\infty}_{+}=\{\xi\in L^{\infty}\,|\,\xi\geqslant 0\}$.

Let us first recall the notion of a random normed module.
\begin{definition}(see \cite{Guo-JFA,Guo-recent})
Let $S$ be a left module over $L^0$ (briefly, an $L^0$-module), a mapping $\|\cdot\|: S\to L^0_+$ is called an $L^0$-norm on $S$ if:

\noindent(i) $\|x\|=0$ if and only if $x=\theta$( the null element of $S$);\\
\noindent(ii)  $\|\xi x\|=|\xi|\|x\|$ for all $\xi\in L^{0}$ and $x\in S$;\\
\noindent(iii) $\|x_1+x_2\|\leqslant \|x_1\|+\|x_2\|$ for all $x_1, x_2\in S$.

\noindent In this case, $(S,\|\cdot\|)$ is called a random normed module (briefly, an RN module).
\end{definition}

In this paper, given an RN module $(S,\|\cdot\|)$, it is always endowed with the $(\varepsilon,\lambda)$-topology. It suffices to know that the $(\varepsilon,\lambda)$-topology is a metrizable linear topology, a sequence $\{x_n\}_{n\geq 1}$ in $S$ converges in the $(\varepsilon,\lambda)$-topology to $x$ iff the sequence $\{\|x_n-x\|\}_{n\geq 1}$ in $L^0_+$ converges in probability to $0$.

We then recall the notion of an $L^\infty$-normed module.

\begin{definition}(\cite{ET1}.)
 Let $E$ be a left module over $L^\infty$ (briefly, an $L^\infty$-module), a mapping $\|\cdot\|: E\to L^\infty_+$ is called an $L^\infty$-norm on $E$ if:

\noindent(i) $\|x\|=0$ if and only if $x=\theta$( the null element of $E$);\\
\noindent(ii)  $\|\xi x\|=|\xi|\|x\|$ for all $\xi\in L^\infty$ and $x\in E$;\\
\noindent(iii) $\|x_1+x_2\|\leqslant \|x_1\|+\|x_2\|$ for all $x_1, x_2\in E$.

\noindent In this case, $(E,\|\cdot\|)$ is called an $L^\infty$--normed module.
\end{definition}

Obviously, if $(E,\|\cdot\|)$ is an $L^\infty$-normed module, then the mapping $\|\cdot\|_\infty: E\to [0,\infty)$ defined by $\|x\|_\infty=\hbox{the essential supremum of~}\|x\|,\forall x\in E$, is a norm on $E$. In this paper we always endow an $L^\infty$--normed module $(E,\|\cdot\|)$ with the topology induced by the norm $\|\cdot\|_\infty$.

\section{The $L^0$-extension of an $L^\infty$-normed module }
 \label{}

Assume that $(S,\|\cdot\|)$ is an RN module, denote $L^\infty(S)=\{x\in S~|~\|x\|\in L^\infty\}$, then it is easy to see that $(L^\infty(S), \|\cdot\|)$ is an $L^\infty$--normed module. For a special case $(S,\|\cdot\|)=(L^0,|\cdot|)$, we have $(L^\infty(S), \|\cdot\|)=(L^\infty,|\cdot|)$. A question is immediately raised: if $(E,\|\cdot\|)$ is an $L^\infty$--normed module, does there exist an RN module $(S,\|\cdot\|)$ such that $(L^\infty(S), \|\cdot\|)=(E,\|\cdot\|)$? We will find the answer in this section.

\subsection{The $L^0$-extension RN module}

Besides the usual norm topology, $L^\infty$ can be also endowed with the topology of convergence in probability. Endowed with this topology, $L^\infty$ is still a metrizable linear topological space, which is typically no longer complete. By the way of completion, we will get $L^0$. This procedure can be done in an abstract way so that we can extend an $L^\infty$-normed module to an RN module.

\begin{theorem}\label{extention}
Let $(E,\|\cdot\|)$ be an $L^\infty$-normed module, then there exist a complete RN module $(E_0, \|\cdot\|_0)$ together with an $L^\infty$-module homomorphism $T: E\to L^\infty(E_0)$ such that: \\
\noindent(1). $\|T(x)\|_0=\|x\|,\forall x\in E$;\\
\noindent(2). $T(E)=\{T(x):x\in E\}$ is a dense subset of $E_0$.
\end{theorem}

{\em proof.}
Define a metric $d$ on $E$ by $d(x,y)=E[1\wedge \|x-y\|], \forall x, y\in E$, where $E[\xi]$ means $\xi$'s expectation with respect to $P$. Obviously, for any sequence $\{x_n\}_{n\geq 1}$ in $E$ and an element $x\in E$, $d(x_n, x)$ tends to $0$ if and only if the sequence $\{\|x_n-x\|\}_{n\geq 1}$ in $L^\infty$ converges to $0$ in probability. It is easy to verify that $(E,d)$ is a metric linear space, namely, the addition operation and scalar multiplication operation are both continuous with respect to the metric $d$. Generally, $(E,d)$ is not complete. Let $E_0$ be the set of all $d$-Cauchy sequences in $E$, where two $d$-Cauchy sequences $\{x_n\}_{n\geq 1}$ and $\{y_n\}_{n\geq 1}$ are identified when $\lim\limits_{n\to\infty}d(x_n,y_n)=0$. As usual, define the addition operation $+:E_0\times E_0\to E_0$ and the scalar multiplication operation $\cdot:L^0\times E_0\to E_0$ by $\{x_n\}_{n\geq 1}+\{y_n\}_{n\geq 1}=\{x_n+y_n\}_{n\geq 1}$ and $a\cdot\{x_n\}_{n\geq 1}=\{ax_n\}_{n\geq 1}$ for any $\{x_n\}_{n\geq 1}, \{y_n\}_{n\geq 1}\in E_0$ and $a\in {\mathbb R}$. Moreover, define a metric $d_0$ on $E_0$ by $d_0(\{x_n\}_{n\geq 1},\{y_n\}_{n\geq 1})=\lim\limits_{n\to\infty}d(x_n,y_n)$. From standard functional analysis, we know that $(E_0,d_0)$ is a complete metric linear space, and $E$ is isometric with a dense subset $\{{\tilde x}: x\in E\}$ of $E_0$, here and in the following $\tilde x$ stands for the constant value sequence $\{x,x,x,\dots\}$ for each $x\in E$.

Now we introduce a module multiplication operation $\ast: L^0\times E_0\to E_0$ and an $L^0$-norm $\|\cdot\|_0: E_0\to L^0_+$ to make $E_0$ become an $L^0$-normed module.

First, we introduce the module multiplication operation. For each $\xi\in L^0$, there exists a sequence $\{\xi_n\}_{n\geq 1}$ in $L^\infty$ such that $\{\xi_n\}_{n\geq 1}$ converges to $\xi$ in probability. For any $d$-Cauchy sequence $\{x_n\}_{n\geq 1}$ in $E$, it is easy to verify that $\{\xi_nx_n\}_{n\geq 1}$ is still a $d$-Cauchy sequence. Moreover, if $\{\eta_n\}_{n\geq 1}$ is another sequence in $L^\infty$ which converges to $\xi$ in probability, and $\{y_n\}_{n\geq 1}$ is another $d$-Cauchy sequence in $E$ such that $\lim\limits_{n\to\infty}d(x_n,y_n)=0$, then we can verify that $\lim\limits_{n\to\infty}d(\xi_nx_n,\eta_ny_n)=0$. Thus, by defining $\xi\ast\{x_n\}_{n\geq 1}=\{\xi_nx_n\}_{n\geq 1}$ for each $\xi\in L^0$ and $\{x_n\}_{n\geq 1}\in E_0$, we get a well-defined mapping $\ast: L^0\times E_0\to E_0$. Further, we can verify that the mapping $\ast$ is indeed a module multiplication operation.

Then, we give the $L^0$-norm. For any $x\in E_0$, there exists a sequence $\{x_n\}_{n\geq 1}$ in $E$ such that $\{{\tilde x}_n\}_{n\geq 1}$ converges to $x$, then $\{\|x_n\|\}_{n\geq 1}$ must converge in probability to some $\eta\in L^0_+$, clearly, this $\eta$ does not depend on the choice of the sequence $\{x_n\}_{n\geq 1}$ in $E$ such that $\{{\tilde x}_n\}_{n\geq 1}$ converges to $x$, namely $\eta$ is uniquely decided by $x$, write $\|x\|_0$ for this unique $\eta$ for each $x$, then we define a mapping $\|\cdot\|_0: E_0\to L^0_+$ by $x\mapsto \|x\|_0$. We can verify that the mapping $\|\cdot\|_0$ is indeed an $L^0$-norm on $E_0$.

Finally, we verify that $d_0(x,y)=E[\|x-y\|_0\wedge 1],~\forall x, y\in E_0$, which means that the $(\varepsilon,\lambda)$-topology induced by $\|\cdot\|_0$ on $E_0$ is the same with that induced by the metric $d_0$. In fact, arbitrarily choose two sequences $\{x_n\}_{n\geq 1}$ and $\{y_n\}_{n\geq 1}$ in $E$ such that $\{{\tilde x}_n\}_{n\geq 1}$ and $\{{\tilde y}_n\}_{n\geq 1}$ converges to $x$ and $y$, respectively, then $d_0(x,y)=\lim\limits_n d(x_n,y_n)=\lim\limits_n E[\|x_n-y_n\|\wedge 1]=E[\lim\limits_n \|x_n-y_n\|\wedge 1]=E[\|x-y\|_0\wedge 1]$. Thus, $(E_0,\|\cdot\|_0)$ is a complete RN module.

To complete the proof, it suffices to define $T:E\to L^\infty(E_0)$ by $Tx=\tilde x, \forall x\in E$.
\hfill\done

\begin{proposition}
Let $(E,\|\cdot\|)$, $(E_0, \|\cdot\|_0)$ and $T: E\to L^\infty(E_0)$ be the same as Theorem \ref{extention} above, then $(E_0, \|\cdot\|_0)$ is unique in the sense of RN module isometric isomorphism, that is to say, if $(S,\|\cdot\|_S)$ is another complete RN module with an $L^\infty$-module homomorphism $R: E\to L^\infty(S)$ such that $\|Rx\|_S=\|x\|,\forall x\in E$ and $R(E)=\{Rx:x\in E\}$ is a dense subset of $S$, then there exists an $L^0$-module isomorphism $U: E_0\to S$ with $\|Ux\|_S=\|x\|_0,\forall x\in E_0$.
\end{proposition}

{\em proof.}
Define $U: E_0\to S$ by $Ux=\lim_{n}Rx_n, \forall x\in E_0$, where $\{x_n\}_{n\geq 1}$ is a sequence in $E$ such that $\{Tx_n\}_{n\geq 1}$ converges to $x$ in $E_0$ and the limit on the right side is taken in $(S,\|\cdot\|_S)$. We can verify that $U$ is an $L^0$-module isomorphism and $\|Ux\|_S=\|x\|_0,\forall x\in E_0$.
\hfill\done

In the sequel, for a given $L^\infty$-normed module $(E,\|\cdot\|)$, we always use $(E_0,\|\cdot\|_0)$ to denote a complete RN module satisfying all the conditions stated in Theorem \ref{extention}, which is in fact unique in the sense of RN module isometric isomorphism according to Proposition 3.2. We call $(E_0,\|\cdot\|_0)$ the $L^0$-extension RN module of $(E,\|\cdot\|)$, and identify $E$ with the dense subset $\{Tx: x\in E\}$ of $E_0$.

Take two examples of $L^\infty$-normed modules given in \cite[Section 8]{ET2}, we give their $L^0$-extension RN modules as follows.

\begin{example}
Let $(\Omega,{\mathcal F},P)$ be a probability space.

(1).Assume that ${\mathcal G}$ is a sub-$\sigma$-algebra of ${\mathcal F}$. For $p\in [1,\infty]$, set $$L^p_{\lambda}({\mathcal F})=\{\xi\in L^0({\mathcal F})~|~\|\xi\|_{p,{\mathcal G}}\in \lambda\}$$
where $\lambda$ stands for $L^\infty({\mathcal G})$ and $\|\xi\|_{p,{\mathcal G}}$ means $\{E[|\xi|^p|{\mathcal G}]\}^\frac{1}{p}$ if $1\leq p<\infty$, or means
$\wedge\{\eta\in {\bar L}^0({\mathcal G}):\eta\geq\xi\}$ if $p=\infty$. Then $(L^p_{\lambda}({\mathcal F}),\|\cdot\|_{p,{\mathcal G}})$ is an $L^\infty({\mathcal G})$-normed module. Its $L^0({\mathcal G})$-extension RN module is $(L^p_{\mathcal G}({\mathcal F}), \|\cdot\|_{p,{\mathcal G}})$ as in \cite[Example 2.5]{FKV}.

(2).Assume that ${\mathcal G}, {\mathcal G}_1$ are two sub-$\sigma$-algebras of ${\mathcal F}$ with ${\mathcal G}\subset {\mathcal G}_1$. Denote $\bar {\mathcal G}=({\mathcal G}, {\mathcal G}_1)$. For ${\bar p}=(p_1,p_2)\in [1,\infty]^2$, set $$L^{\bar p}_\lambda(\bar {\mathcal G})=\{\xi\in L^0({\mathcal F})~|~\|\xi\|_{{\bar p},\bar {\mathcal G}}:=\left\|\|\xi\|_{p_2,{\mathcal G}_1}\right\|_{p_1,{\mathcal G}}\in \lambda\}$$
where $\lambda$ stands for $L^\infty({\mathcal G})$.
Then $(L^{\bar p}_\lambda(\bar {\mathcal G}),\|\cdot\|_{{\bar p},\bar {\mathcal G}})$ is an $L^\infty({\mathcal G})$-normed module. Its $L^0({\mathcal G})$-extension RN module is $(L^{\bar p}_{\bar {\mathcal G}}(\mathcal F), \|\cdot\|_{{\bar p},\bar {\mathcal G}})$, where $L^{\bar p}_{\bar {\mathcal G}}(\mathcal F):=\{\xi\in L^0({\mathcal F})~|~\|\xi\|_{{\bar p},\bar {\mathcal G}}:=\left\|\|\xi\|_{p_2,{\mathcal G}_1}\right\|_{p_1,{\mathcal G}}\in L^0(\mathcal G)\}$.

\end{example}

\subsection{Two sub-$L^0$-modules of $E_0$}

We remind that we identify $E$ with the dense subset $\{Tx: x\in E\}$ of $E_0$.

Let $L(E)=L^0\ast E:=\{\xi\ast x: \xi\in L^0, x\in E\}$, since $\xi\ast(\eta\ast x)=(\xi\eta)\ast x$ and $\xi\ast x+\eta\ast y=\gamma\ast(\frac{\xi}{\gamma}x+\frac{\eta}{\gamma}y)$, where $\gamma=|\xi|+|\eta|+1\in L^0_{++}$, we see that $L(E)$ is in fact the smallest sub-$L^0$-module in $E_0$ which contains $E$.

 For the brevity, we introduce a notation. In the sequel, denoted by $$\Pi({\mathcal F})=\{\{A_n: n\in \N\}: A_n\in {\mathcal F}\mbox{~for all~} n;~ \cup_n A_n=\Omega;~ A_n\cap A_m=\emptyset \mbox{~for all~}m\neq n\}$$ the set of all countable partitions of $\Omega$ to ${\mathcal F}$.

For each sequence $\{x_n,n\in \N\}$ in $E_0$ and a countable partition $\{A_n: n\in \N\}\in \Pi({\mathcal F})$, there always exists one and only one $x$ in $E_0$ such that ${\tilde I}_{A_n}x={\tilde I}_{A_n}x_n, \forall n\in \N$. For the existence, take $x^\prime_k=\sum_{i=1}^k{\tilde I}_{A_i}x_i$ for each $k\in \N$, then $\{x^\prime_k,k\in \N\}$ forms a Cauchy sequence in $E_0$, thus converges to some $x$ in $E_0$, further for this $x$, ${\tilde I}_{A_n}x=\lim_k{\tilde I}_{A_n}x^\prime_k={\tilde I}_{A_n}x_n, \forall n\in \N$. For the uniqueness, assume that $x^\prime$ also satisfies that ${\tilde I}_{A_n}x={\tilde I}_{A_n}x_n, \forall n\in \N$, then $x-x^\prime=\lim_k{\tilde I}_{\cup^k_{i=1} A_i}(x-x^\prime)=0$. We write $\sum_{n\in N}{\tilde I}_{A_n}x_n$ for the unique $x$ such that ${\tilde I}_{A_n}x={\tilde I}_{A_n}x_n, \forall n\in \N$.

We recall the notion of the countable concatenation property introduced by Guo\cite{Guo-JFA} as follows.

\begin{definition}\label{ccp}
Assume that $G$ is a subset of an RN module, we say that $G$ has the countable concatenation property, if there exists $g\in G$ for each sequence $\{g_n:n\in \N\}$ in $G$ and countable partition $\{A_n: n\in \N\}\in \Pi({\mathcal F})$ such that ${\tilde I}_{A_n}g={\tilde I}_{A_n}g_n,~\forall n\in \N$.
\end{definition}

From the argument just before Definition \ref{ccp}, $E_0$ has the countable concatenation property. Given a subset $G$ of $E_0$, denote $$H^0_{cc}(G)=\{\sum_{n\in \N}{\tilde I}_{A_n}g_n: g_n \in G~\mbox{for all}~ n;~ \{A_n: n\in \N\}\in \Pi({\mathcal F})\},$$
then $G$ has the countable concatenation property if and only if $H^0_{cc}(G)=G$.

\begin{proposition}
$H^0_{cc}(E)$ is an sub-$L^0$-module of $E_0$.
\end{proposition}
{\em proof.}
For each $x=\sum_{i\in \N}{\tilde I}_{A_i}x_i$ and $y=\sum_{j\in \N}{\tilde I}_{B_j}y_j$ in $H^0_{cc}(E)$, where $\{x_i, i\in \N\}$ and $\{y_j,j\in \N\}$ in $E$, and $\{A_i: i\in \N\}$ and $\{B_j: j\in \N\}$ in $\Pi({\mathcal F})$, then $\{A_i\cap B_j:i,j\in \N\}\in \Pi({\mathcal F})$ and ${\tilde I}_{A_i\cap B_j}(x+y)={\tilde I}_{A_i\cap B_j}({\tilde I}_{A_i}x+{\tilde I}_{B_j}y)={\tilde I}_{A_i\cap B_j}({\tilde I}_{A_i}x_i+{\tilde I}_{B_j}y_j)\in E$ for each $i,j\in \N$, which means that $x+y\in H^0_{cc}(E)$.

For each $\xi\in L^0$ and $x=\sum_{i\in \N}{\tilde I}_{A_i}x_i\in H^0_{cc}(E)$, where $\{x_i: i\in \N\}\subset E$ and $\{A_n: n\in \N\}\in \Pi({\mathcal F})$, let $B_j=\{\omega\in\Omega: j-1\leq |\xi^0(\omega)|<j\}$ for each $j\in \N$, where $\xi^0$ is any representative for $\xi$, then $\{B_j:j\in \N\}\in \Pi({\mathcal F})$ and $\{A_i\cap B_j:i,j\in \N\}\in \Pi({\mathcal F})$, then ${\tilde I}_{A_i\cap B_j}\xi\ast x=({\tilde I}_{B_j}\xi) ({\tilde I}_{A_i}x)=({\tilde I}_{A_i\cap B_j}\xi) x_i\in E$ for each $i,j\in \N$, which means that $\xi\ast x\in H^0_{cc}(E)$.
\hfill\done

Since $E\subset H^0_{cc}(E)$ and $L(E)$ is the smallest sub-$L^0$-module in $E_0$ which contains $E$, we obtain that $L(E)\subset H^0_{cc}(E)$.

In the following, we give two examples to show that both the two inclusions $L(E)\subset H^0_{cc}(E)$ and $H^0_{cc}(E)\subset E_0$ may be proper.

\begin{example}\label{Lproper}
Let $\Omega=\{1,2,3,\dots\}$, ${\mathcal F}$ the $\sigma$--algebra consisting of all subsets of $\Omega$, and probability $P$ defined by $P(\{k\})=\frac{1}{2^k}$ for each $k\in \Omega$. Take $$E=\{\phi:\Omega\to {\mathbb R}:~\mbox{there exist some $n_\phi\in \N$ such that}~\phi(k)=0,\forall~ k\geq n_{\phi}\},$$ then $E\subset L^\infty$, and $(E,|\cdot|)$ is an $L^\infty$-normed module. For any $x\in L^0$, take $\phi_n={\tilde I}_{\{1,2,\dots,n\}}x$ for each $n\in \N$, then it is obviously that $\phi_n$ converges to $x$ in probability, so $E_0=L^0$. It is easy to see that $L(E)=E$ is a proper subset of $H^0_{cc}(E)=L^0$.
\end{example}

\begin{example}
Let $\Omega=[0,1], {\mathcal F}={\mathcal B}[0,1], P=\mbox{Lebesgue measure}~\lambda$. $E=\{(\xi_1,\xi_2,\dots,\xi_n,0,0,\dots):n=1,2,\dots;~\xi_i\in L^\infty, i=1,2,\dots,n\}$, $L^\infty$--norm on $E$ is defined by $\|(\xi_1,\xi_2,\dots,\xi_n,0,0,\dots)\|=\vee\{|\xi_k|: k=1,2,\dots,n\}$, then $\{x^n:=({\mathds 1}, \frac{\mathds 1}{2},\dots,\frac{\mathds 1}{n},0,0,\dots)\}_{n\geq 1}$ forms an Cauchy sequence in $E$, clearly, it ``converges'' to $x=({\mathds 1}, \frac{\mathds 1}{2},\dots,\frac{\mathds 1}{n},\dots)$, which implies $x\in E_0$, however $x\notin H^0_{cc}(E)$ since $\tilde I_Ax$ is not in $E$ for any $A\in {\mathcal F}$. Thus $H^0_{cc}(E)$ is a proper subset of $E_0$.
\end{example}

\begin{proposition}\label{Equa}
Assume that $(E,\|\cdot\|)$ is complete, then $E_0=H^0_{cc}(E)=L(E)$.
\end{proposition}

{\em proof.}
For any fixed $x\in E_0$, we will express $x$ as a countable concatenation. According to the definition of $E_0$, there exists a sequence $\{x_n\}_{n\geq 1}$ in $E$ such that the sequence $\{\|x_n-x\|\}_{n\geq 1}$ in $L^0_+$ converges to $0$ in probability.
By passing to an appropriate subsequence if necessary, we can assume that the sequence $\{\|x_n-x\|\}_{n\geq 1}$ converges to 0 almost surely. Then according to Egoroff's theorem, we can find $A_1\in {\mathcal F}$ with $P(A_1)\geq \frac{1}{2}$ such that $\{{\tilde I}_{A_1}\|x_n-x\|\}_{n\geq 1}$ converges to $0$ uniformly, using Egoroff's theorem once again, we can find $A_2\in {\mathcal F}$ with $A_2\subset \Omega\setminus A_1,~P(A_2)\geq \frac{1}{2}(1-P(A_1))$ such that $\{{\tilde I}_{A_2}\|x_n-x\|\}_{n\geq 1}$ converges to $0$ uniformly. By induction, if $A_1\cup A_2\cup\cdots\cup A_k\neq \Omega$, then by using Egoroff's theorem we can obtain $A_{k+1}\in {\mathcal F}$ with $A_{k+1}\subset \Omega\setminus(A_1\cup A_2\cup\cdots\cup A_k), P(A_{k+1})\geq \frac{1}{2}(1-P(A_1)-\cdots-P(A_k))$ such that $\{{\tilde I}_{A_{k+1}}\|x_n-x\|\}_{n\geq 1}$ converges to $0$ uniformly. To sum up,  we get a countable partition $\{A_k:k\in \N\}\in \Pi({\mathcal F})$ such that $\{{\tilde I}_{A_k}\|x_n-x\|\}_{n\geq 1}$ converges to $0$ uniformly for every $A_k$. For each $k$, it is clear that $\{{\tilde I}_{A_k}x_n\}_{n\geq 1}$ is a Cauchy sequence in $E$, thus converges to some $y_k\in E$ due to the assumption that $E$ is complete. Then we have $x=\sum^\infty_{k=1}{\tilde I}_{A_k}y_k$. Thus $E_0=H^0_{cc}(E)$.

We then express $x=\sum^\infty_{k=1}{\tilde I}_{A_k}y_k$ as $x=\xi \ast y$ for some $\xi\in L^0$ and $y\in E$. Let $$y^\prime_k=\sum^k_{i=1}{\tilde I}_{A_i}\frac{1}{(1+\|y_i\|_\infty)2^i}y_i$$
Obviously, $\{y^\prime_k\}_{k\geq 1}$ forms a Cauchy sequence in $E$, thus converges to some $y\in E$ since $E$ is complete. Set $\xi=\sum^\infty_{i=1}{(1+\|y_i\|_\infty)2^i}{\tilde I}_{A_i}\in L^0$, it is easy to verify that $x=\sum^\infty_{k=1}{\tilde I}_{A_k}y_k=\xi\ast y$.
\hfill\done

\subsection{Recover $E$ from $E_0$}

 In Theorem \ref{extention}, we always have an inclusion $E\subset L^\infty(E_0)$. Does the inverse inclusion $L^\infty(E_0)\subset E$ also hold? equivalently, do we have $E= L^\infty(E_0)$?

Let us return to Example \ref{Lproper}. We have $E_0=L^0$, then $L^\infty(E_0)=L^\infty$, thus $E$ is merely a proper subset of $L^\infty(E_0)$. Therefore, the answer to the above question is negative. We then look for conditions to make the equality $E=L^\infty(E_0)$ hold.

Let us have some observations. Obviously, $(L^\infty(E_0),\|\cdot\|)$ is always a complete $L^\infty$-normed module. Moreover, we can see that the unit ball $U(L^\infty(E_0)):=\{x\in L^\infty(E_0): \|x\|\leqslant 1\}$ always has the countable concatenation property. Thus, if $E=L^\infty(E_0)$ holds, then $E$ must be complete and its unitary ball $U(E):=\{x\in E: \|x\|\leqslant 1\}$ must have the countable concatenation property. The following theorem says that these two necessary conditions are also sufficient.

\begin{proposition}\label{Equality}
Assume that $(E,\|\cdot\|)$ is an complete $L^\infty$-normed module, and its unitary ball $U(E)$ has the countable concatenation property, then $E=L^\infty(E_0)$.
\end{proposition}
{\em proof.}
Since $(E,\|\cdot\|)$ is complete, it follows from Proposition \ref{Equa} that $E_0=H^0_{cc}(E)$. For any $x\in L^\infty(E_0)$, there exist $\{A_n: n\in \N\}\in \Pi({\mathcal F})$ and $\{x_n,n\in \N\}$ in $E$, such that ${\tilde I}_{A_n}x={\tilde I}_{A_n}x_n, \forall n\in \N$. Set $\lambda =\|x\|_\infty+1$, then for each $n$, we have $\|{\tilde I}_{A_n}\lambda ^{-1}x_n\|=\|{\tilde I}_{A_n}\lambda ^{-1}x\|\leqslant 1$, namely, ${\tilde I}_{A_n}\lambda ^{-1}x={\tilde I}_{A_n}\lambda ^{-1}x_n\in U(E)$. Using the assumption that $U(E)$ has the countable concatenation property, we get $\lambda ^{-1}x\in U(E)$, thus $x\in E$.
\hfill\done

\subsection{The dual $L^\infty$-normed module of $E$ and the random conjugate space of $E_0$}

Let $(E,\|\cdot\|)$ be an $L^\infty$-normed module, denote by $E^\prime$ the $L^\infty$-module of all continuous module homomorphisms $f:(E,\|\cdot\|)\to (L^\infty,|\cdot|)$, and define $\|\cdot\|^\prime: E^\prime\to L^\infty_+$ by $\|f\|^\prime=\vee\{|f(x)|:x\in E, \|x\|\leqslant 1\}$ for each $f\in E^\prime$, then according to \cite[Proposition 9.1]{ET2}, $(E^\prime,\|\cdot\|^\prime)$ is a complete $L^\infty$-normed module, called the dual $L^\infty$-normed module of $(E,\|\cdot\|)$.

Let $(S,\|\cdot\|)$ be an RN module, denote by $S^\ast$ the $L^0$-module of all continuous module homomorphisms $f:(S,\|\cdot\|)\to (L^0,|\cdot|)$, and define $\|\cdot\|^\ast: S^\ast\to L^0_+$ by $\|f\|^\ast=\vee\{|f(x)|:x\in S, \|x\|\leqslant 1\}$ for each $f\in S^\ast$, then $(S^\ast,\|\cdot\|^\ast)$ is a complete RN module, called the random conjugate space of $(S,\|\cdot\|)$.

Assume that $(E,\|\cdot\|)$ is an $L^\infty$-normed module, and $(E^\prime, \|\cdot\|^\prime)$ is its dual $L^\infty$--normed module. Let $(E_0, \|\cdot\|_0)$ be the $L^0$-extension RN module of $(E,\|\cdot\|)$. Then what is the $L^0$-extension RN module of $(E^\prime, \|\cdot\|^\prime)$? The following theorem says that the $L^0$-extension RN module of $(E^\prime, \|\cdot\|^\prime)$ is exactly $((E_0)^\ast,\|\cdot\|^\ast_0)$, the random conjugate space of $(E_0, \|\cdot\|_0)$.

\begin{theorem}\label{Homo}
There exists an $L^\infty$-module isomorphism $U: (E^\prime, \|\cdot\|^\prime)\to (L^\infty (E_0^\ast), \|\cdot\|^\ast_0)$ with $\|Uf\|^\ast_0=\|f\|^\ast_0,~\forall f\in E^\prime$. Thus the $L^0$-extension RN module of $(E^\prime, \|\cdot\|^\prime)$ is $((E_0)^\ast,\|\cdot\|^\ast_0)$.
\end{theorem}
{\em proof.} For each $f\in E^\prime$, define $Uf: (E_0,\|\cdot\|_0) \to (L^0,|\cdot|)$ by $Uf(x)=\lim_n f(x_n)$ for every $x\in E_0$, where $\{x_n\}_{n\geq 1}$ is a sequence in $E$ which converges to $x$ in $E_0$ and the limit on the right side is taken with respect to the topology of convergent in probability, then we can check that $Uf\in E^\ast_0$ and $\|Uf\|^\ast_0=\|f\|^\prime$, thus $Uf\in L^\infty (E_0^\ast)$. It is easy to see that the mapping $U: E^\prime\to L^\infty (E_0^\ast), f\mapsto Uf$ is an $L^\infty$-module homomorphism.

It remains to show that $U$ is surjective. To this end, for each $g\in L^\infty (E_0^\ast)$, let $Rg$ be the restriction of $g$ to $E\subset E_0$, we get $Rg\in E^\prime$ and $\|Rg\|^\prime=\|g\|^\ast_0$. Clearly we have $U(Rg)=g$.
\hfill\done

\section{The reflexivity and subreflexivity of an $L^\infty$-normed module}

For the brevity, for any given $L^\infty$-normed module or RN module $(X,\|\cdot\|)$, we use $U(X):=\{x\in X~|~\|x\|\leqslant 1\}$ to denote the unit ball of $X$.

We recall that: a complete $L^\infty$-normed module $(E,\|\cdot\|)$ is said to be reflexive \cite{ET1}, if $j(U(E))=U(E^{\prime\prime})$, where $j: E\to E^{\prime\prime}$ is the canonical embedding defined by for each $x\in E$, $[j(x)](f)=f(x),\forall f\in E^\prime$; and $(E,\|\cdot\|)$ is said to be subreflexive, if $A(E)=\{f\in E^\prime~|~\mbox{there exists $x\in U(E)$ such that~}|f(x)|=\|f\|^\prime\}$ is dense in $(E^\prime,\|\cdot\|^\prime)$.

Accordingly, a complete RN module $(S,\|\cdot\|)$ is said to be random reflexive, if $j_0(U(S))=U(S^{\ast\ast})$, where $j_0: S\to S^{\ast\ast}$ is the canonical embedding defined by for each $x\in S$, $[j_0(x)](f)=f(x),\forall f\in S^\ast$; and $(S,\|\cdot\|)$ is said to be random subreflexive \cite{Zhao-Guo}, if $A(S)=\{f\in S^\ast~|~\mbox{there exists $x\in U(S)$ such that~}|f(x)|=\|f\|^\ast\}$ is dense in $(S^\ast,\|\cdot\|^\ast)$.

\subsection{The reflexivity of an $L^\infty$-normed module}

For any given sequence $\{f_n\}_{n\geq 1}$ in $U(E^\prime)$ and $\{A_n:n\in \N\}\in \Pi({\mathcal F})$, $\sum^\infty_{n=1}{\tilde I}_{A_n}f_n$ is still an element in $U(E^\prime)$, where $\sum^\infty_{n=1}{\tilde I}_{A_n}f_n: E\to L^\infty$ is defined by $(\sum^\infty_{n=1}{\tilde I}_{A_n}f_n)(x)=\sum^\infty_{n=1}{\tilde I}_{A_n}[f_n(x)]$ for each $x\in E$, which means that $U(E^\prime)$ always has the countable concatenation property. Thus, if $E$ is reflexive, $U(E)$ must have the countable concatenation property, since $U(E^{\prime\prime})$ has. This simple observation gives a necessary condition for a complete $L^\infty$-normed module to be reflexive.

\begin{proposition}\label{connection}
Let $(E,\|\cdot\|)$ be a complete $L^\infty$--module with $U(E)$ having the countable concatenation property and $(E_0, \|\cdot\|_0)$ its $L^0$-extension RN module. Then $E$ is reflexive if and only if $E_0$ is random reflexive.
\end{proposition}
{\em proof.}
According to Proposition \ref{Equality}, $U(E)=U(E_0)$, $U(E^\prime)=U((E^\prime)_0)$ and $U(E^{\prime\prime})=U((E^{\prime\prime})_0)$. By Theorem \ref{Homo}, $((E_0)^\ast,\|\cdot\|^\ast_0)=((E^\prime)_0, \|\cdot\|^\prime_0)$, thus $jU(E)=U(E^{\prime\prime})$ iff $j_0U(E_0)=U((E_0)^{\ast\ast})$, equivalently, $E$ is reflexive iff $E_0$ is random reflexive. \hfill\done

\begin{remark}
In functional analysis, it is well-known that a Banach space $B$ is reflexive iff its dual space $B^\prime$ is reflexive, however, for a complete $L^\infty$-normed module $E$, it is possible that the dual $L^\infty$-normed module $E^\prime$ is reflexive while $E$ is not, considering Example 3.6 for instance. From this point of view, according to Proposition 4.1, the condition that the unit ball of an $L^\infty$--module has the countable concatenation property is another kind of ``completeness'' to some extent.
\end{remark}

\begin{theorem}\label{James}
Let $(E,\|\cdot\|)$ be a complete $L^\infty$-normed module and $(E^\prime,\|\cdot\|^\prime)$ its dual $L^\infty$-normed module. If $U(E)$ has the countable concatenation property, then $E$ is reflexive if and only if there exists an $x_f\in U(E)$ for each $f\in E^\prime$ such that $f(x_f)=\|f\|^\prime$.
\end{theorem}
{\em proof.}
By Proposition \ref{connection}, we have that $E$ is reflexive if and only if $E_0$ is random reflexive, meanwhile, by James' theorem for complete $L^0$-normed modules\cite[Theorem 3.1]{Guo-Li}, $E_0$ is random reflexive if and only if there exists an $x_g\in U(E_0)$ for each $g\in (E_0)^\ast$ such that $g(x_g)=\|g\|^\ast_0$. Thus, it is equivalent to show ``there exists an $x_g\in U(E_0)$ for each $g\in (E_0)^\ast$ such that $g(x_g)=\|g\|^\ast_0$'' if and only if ``there exists an $x_f\in U(E)$ for each $f\in E^\prime$ such that $f(x_f)=\|f\|^\prime$''.

Due to the assumption that $U(E)$ has the countable concatenation property, Proposition \ref{Equality} yields that $U(E)=U(E_0)$, thus if there exists an $x_f\in U(E_0)$ for each $f\in (E_0)^\ast$ such that $f(x_f)=\|f\|^\ast_0$, then for any given $f\in E^\prime\subset (E_0)^\ast$, one can find an $x_f\in U(E_0)=U(E)$ such that $f(x_f)=\|f\|^\ast_0=\|f\|^\prime$. Conversely, if there exists an $x_f\in U(E)$ for each $f\in E^\prime$ such that $f(x_f)=\|f\|^\prime$, since for any given $g\in (E_0)^\ast$, there exists $\{A_n:n\in \N\}\in \Pi({\mathcal F})$ such that ${\tilde I}_{A_n}g\in L^\infty ((E_0)^\ast)=E^\prime$, then for each $n$, there exists $x_n\in U(E)$ such that  ${\tilde I}_{A_n}g(x_n)=\|{\tilde I}_{A_n}g\|^\prime$, thus if we take $x_g=\sum^\infty_{n=1}{\tilde I}_{A_n}x_n\in U(E)=U(E_0)$, then we have $g(x_g)=\|g\|^\ast_0$.
\hfill\done

\begin{corollary}
Let $(E,\|\cdot\|)$ be a complete $L^\infty$-normed module and $(E_0, \|\cdot\|_0)$ its $L^0$-extension RN module. Then if there exists an $x_f\in U(E)$ for each $f\in E^\prime$ such that $f(x_f)=\|f\|^\prime$, $(E_0, \|\cdot\|_0)$ must be random reflexive.
\end{corollary}
{\em proof}
Consider the complete $L^\infty$-normed module $(L^\infty(E_0), \|\cdot\|_0)$, its unit ball has the countable concatenation property and for each $f\in (L^\infty(E_0))^\prime=E^\prime$, we trivially have $x_f\in U(E)\subset U(E_0)=U(L^\infty(E_0))$ such that $f(x_f)=\|f\|^\prime$, thus $(L^\infty(E_0), \|\cdot\|_0)$ is reflexive by Theorem \ref{James}, then $(E_0, \|\cdot\|_0)$ must be random reflexive follows from Proposition \ref{connection}.
\hfill\done

\begin{remark}
For a complete $L^\infty$-normed module $(E,\|\cdot\|)$, we are not clear whether the statement ``there exists an $x_f\in U(E)$ for each $f\in E^\prime$ such that $f(x_f)=\|f\|^\prime$'' implies that ``$(E,\|\cdot\|)$ is reflexive'', or implies that ``U(E) has the countable concatenation property'', namely, we are not sure whether James theorem holds in complete $L^\infty$-normed modules or not.
\end{remark}

\subsection{The subreflexivity of an $L^\infty$-normed module}

\begin{theorem}
 Let $(E,\|\cdot\|)$ be a complete $L^\infty$-normed module with $U(E)$ having the countable concatenation property, then $E$ is subreflexive.
\end{theorem}
{\em proof.}
Obviously, $U(E)=U(E_0)$ is a closed $L^0$-convex and a.s bounded (by definition, $\|x\|_0\leqslant 1$ for all $x\in U(E)$) subset of $E_0$, immediately, $U(E)$ is a bounded, closed and convex subset of the Banach space $(L^1(E_0), \|\cdot\|_1)$, thus $$A_1(E_0)=\{f\in (L^1(E_0))^\prime~|~\mbox{there exists an}~x\in U(E)~\mbox{such that}~f(x)=\|f\|\}$$ is dense in $L^1(E_0))^\prime$ by the Bishop-Phelps theorem, using the fact that $(L^1(E_0))^\prime\cong L^\infty(E^\ast_0)=E^\prime$ and \cite[Lemma 3.3]{Wu1}, we conclude that $A(E)\cong A_1(E_0)$ is dense in $E^\prime$.
\hfill\done

We give a example to show that, without the assumption that $U(E)$ has the countable concatenation property, a complete $L^\infty$-normed module $(E,\|\cdot\|)$ may be not subreflexive.

\begin{example}
Let $\Omega=\{1,2,3,\dots\}$, ${\mathcal F}$ the $\sigma$-algebra which consists of all subsets of $\Omega$, probability $P$ defined by $P(\{k\})=\frac{1}{2^k}$ for each $k\in \Omega$. Let $E=\{\phi:\Omega\to R: \lim\limits_{k\to \infty}\phi(k)=0\}$, then $E$ is a closed sub-$L^\infty$-module of $(L^\infty, |\cdot|)$, thus $(E, |\cdot|)$ is a complete $L^\infty$-normed module. It is easy to verify that $E^\prime=L^\infty$ and $A(E)=E$. Since $E$ is a proper closed sub-$L^\infty$-module of $L^\infty$, we conclude that $E$ is not subreflexive.
\end{example}





\end{document}